\documentclass[10pt]{amsart}
\usepackage[cp1251]{inputenc}
\usepackage[english,russian]{babel}
\usepackage{amsmath}
\usepackage{amssymb}
\usepackage{amsfonts}
\usepackage{graphicx}

\newtheorem{lemma}{Lemma}
\newtheorem{theorem}{Theorem}

\newtheorem{proposition}{Proposition}

\begin{document}
\renewcommand{\refname}{References}
\renewcommand{\proofname}{Proof.}
\renewcommand{\figurename}{Fig.}

\thispagestyle{empty}

\title[Modifications of Simon text model]{Modifications of Simon text model}
\author{{M.G. Chebunin, A.P. Kovalevskii}}%

\thanks{\sc Chebunin, M.G., Kovalevskii, A.P.,
Modifications of Simon text model}
\thanks{\rm The reported study was funded by RFBR and CNRS according to the research project No. 19-51-15001.}

\maketitle {\small
\begin{quote}
\noindent{\sc Abstract. } 
We discuss probability text models and their modifications.
We construct processes of different and unique words in a text.
The models are to correspond to the real text statistics.
The infinite urn model (Karlin model) and the Simon model are the most known models of texts, but they do not give the ability to simulate the number of unique words correctly. The infinite urn model give sometimes the incorrect limit of the relative number of unique and different words. The Simon model states a linear growth of the numbers of different and unique words.
We propose three modifications of the Karlin and Simon models. The first one is the offline variant, the Simon model starts after the completion of the infinite urn scheme. We prove limit theorems for this modification in embedded times only. The second modification involves the compound Poisson process in the infinite urn model. We prove limit theorems for it.
The third modification is the online variant, the Simon redistribution works at any toss of the Karlin model. In contrast to the compound Poisson model, we have no analytics for this modification.
We test all the modifications by the simulation and have a good correspondence to the real texts.
\medskip

 \end{quote}
}

\section{Introduction}

Probabilistic text modeling involves several simplifications. However, the probabilistic model should maintain the behavior of text statistics that are observed in practice.
In particular, we will consider the number of different words in the text and the number of words that occur once.

Let $R_n$ be a number of different words among the first $n$ words of the text.
$R_{n,i}$ be a number of words encountered $i$ times,
$R^*_{n,i}$ be the number of words encountered not lesser than $i$ times.
Therefore $R_n=R^*_{n,1}$, $R_{n,i}=R^*_{n,i}-R^*_{n,i+1}$, $i\ge 1$.

The power law of the growth of the number of different words is called Herdan's Law or Heaps' Law. 
It refers to Herdan \cite{18} and Heaps \cite{17}.

Bahadur \cite{3} and Karlin \cite{21} studied an infinite urn model: any new ball goes to some of infinitely many urns
with probability that corresponds to a power law and independently of anything else. 
The interpretation for texts corresponds
words of a text to balls and words of a dictionary to urns.

Simon \cite{25} proposed quite another model: the $(n+1)$-th word of a text is a new one with probability $p$,
 and coinсides with any of previous words with probability $(1-p)/n$.

The infinite urn scheme looks more suitable for describing real texts, since Simon's model
leads to a linear increase in the number of different words. However, the infinite urn scheme is not flexible enough
to describe texts. We study two different estimates for the parameter $ \theta $ of the exponential decay of the probabilities.
One of them is $ \widehat{\theta} $, it characterizes the rate at which the number of different words grows. Another estimate $ \theta^{*} $
is the ratio of the number of unique words to the number of different words. According to the infinite urn scheme with exponential decay of probabilities, these two estimates should converge to the same number $ \theta $.

But we show in Section 2 the examples that the estimates are substantially different, the number of words encountered once (unique words)
grows according to a power law with the same exponent but with a lower constant.

So we need some modifications or combinations of the models.

In Section 3, we study an elementary probabilistic text model. This is an infinite urn  scheme.  In this model, the number of
different words and the number of words encountered once  were studied by Bahadur\cite{3}, Karlin \cite{21}, 
Chebunin and Kovalevskii \cite{10}. We study the correspondence between the empirical and theoretical behavior of these statistics.

In Section 4, we study Simon's model.
Simon \cite{25}  proposed the next stochastic model: the $ (n + 1) $-th word in the text is new with probability  $ p $; 
it coincides with each of the previous words with probability
$(1-p)/n$. In fact, he proposed a more general model with the same dynamics of numbers of word occurrences. 
He based his model on the model of Yule \cite{26} who constructed it to explain the distribution of
biological genera by number of species. 
 Baur and Bertoin \cite{8} proposed a modification of the Yule-Simon model with a wide class of limiting distributions.

We study the asymptotic behavior of the statistics $ R_{n, 1} $ in the Simon model based on
functional limit theorems for urn models obtained by Janson.

In Section 5, we propose the offline Simon modification of the Karlin model. The purpose of these modifications is to correspond
the theoretical and empirical behavior of the sequences $ \{R_j \}_{j \ge 1} $ and $ \{R_{j, 1} \}_{j \ge 1} $. 
 We prove analytical theorems (SLLN and FCLT) for the process in embedded times of increation of the initial urn process.

In Section 6, we study the second modification. It involves the compound Poisson process in the infinite urn model.  We prove SLLN and FCLT for the modification. 
In Section 7, we propose the third modification. It is the online variant, the Simon redistribution works at any toss of the Karlin model.  
In contrast to the compound Poisson model, we have no limit theorems for this modification.

We test all the modifications by the simulation and have a good correspondence to the real texts.
We discuss the advantages and disadvantages of the models in Section 8.

\begin{figure}[htbp]
\begin{center}
\includegraphics[bb= 0 0 801 240]{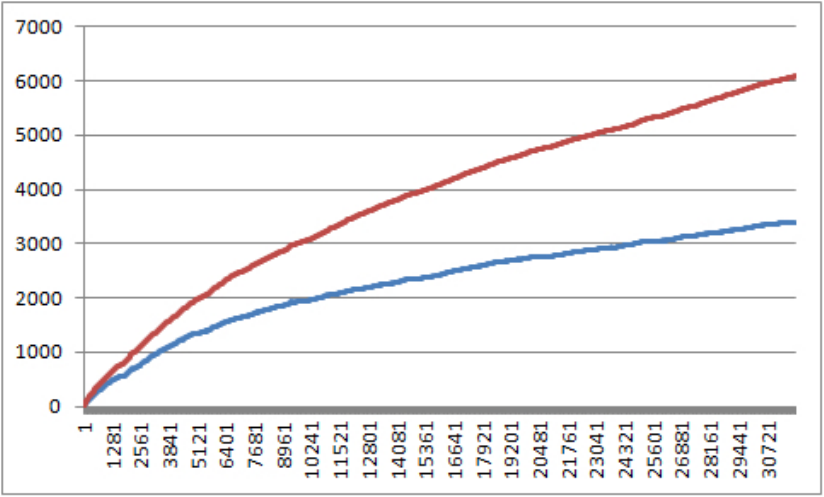}
\caption{ \centering Numbers of different words (red) and numbers of unique words (blue) in {\it Childe Harold's Pilgrimage}  
by Byron }
\end{center}
\end{figure}

\section{Empirical analysis}

We analyze the number of different words and the number of words encountered once in texts of different authors.
These two processes behave like power functions with the same exponent
but with different factors.

We estimate the exponent $ \theta \in (0,1) $ of the power functions in two different ways.
Chebunin and Kovalevskii \cite{12} proposed the estimate 
\[
\widehat{\theta}=\log_2 R_n - \log_2 R_{[n/2]}
\]
and studied conditions of its consistency and asymptotic normality.
It has been used for the analysis of short texts (Zakrevskaya and Kovalevskii \cite{27}). 

Another estimate is 
\[
\theta^{*}=R_{n,1}/R_n.
\]

Karlin \cite{21} proved that it is consistent for the elementary text model under weak assumptions (see the next section).   
This is the asymptotically normal estimate under some additional assumptions (Chebunin and Kovalevskii \cite{11}).

For {\it Childe Harold's Pilgrimage}  
by Byron we have $n=37064$, $R_n=6911$, $R_{n/2}=4582$, $\widehat{\theta}=0.5929$, $R_{n,1}=3912$, 
$\theta^{*}=0.5661$, so the second estimate is significantly smaller than the first one.

\begin{figure}[htbp]
\begin{center}
\includegraphics[bb= 0 0 801 230]{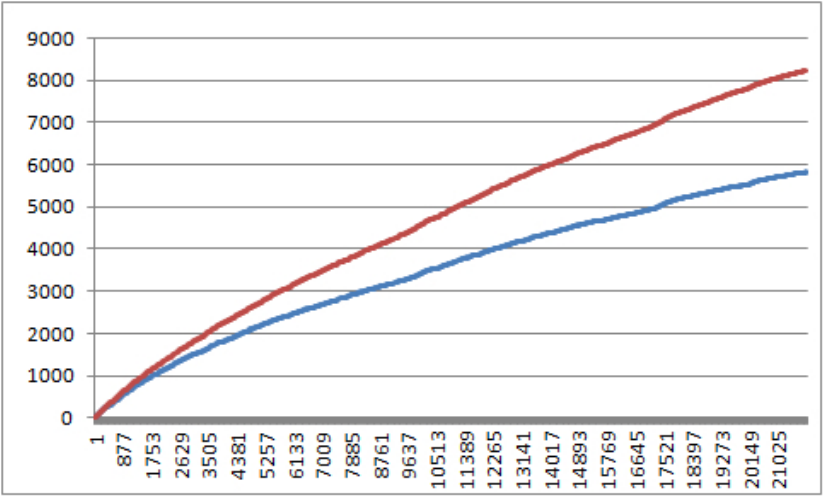}
\caption{ \centering Numbers of different words (red) and numbers of unique words (blue) in {\it Evgene Onegin}  
by Pushkin}
\end{center}
\end{figure}

For {\it Evgene Onegin}  
by Pushkin (in Russian) we have $n=21882$, $R_n=8236$, $R_{n/2}=4916$, $\widehat{\theta}=0.7445$, $R_{n,1}=5824$, 
$\theta^{*}=0.7071$, so the second estimate is significantly smaller than the first one.

\section{Results for the elementary urn model}

The simplest probabilistic model of text is the infinite urn scheme.
Words are selected sequentially independently of each other from an infinite dictionary.
The probabilities of the appearance of words decrease in accordance with the power distribution according to Zipf's law.
As Bahadur showed, the number of different words is growing according to a power law. Karlin showed that the number of words
met once grows under this model also according to a power law with the same exponent.

Karlin \cite{21} studied an infinite urn scheme, that is, $n$ balls distributed to urns independently and randomly;
there are infinitely many urns. Each ball goes to urn $i$ with probability $p_i>0$, $p_1+p_2+\ldots=1$
(without loss of generality $p_1\geq p_2 \geq \ldots $). 

Let (see Karlin \cite{21}) $ \Pi =\{\Pi(t),\ t \geq 0\} $ be a Poisson process with parameter 1.
We denote by $ X_i (n) $ a number of balls in urn $ i $.
According to well-known property of splitting of Poisson flows,
stochastic processes
$\{X_i(\Pi(t)), \ t\geq 0\}$  are Poisson with
 intensities $ p_i $ and are mutually independent for different $ i $'s. The definition implies that
\[
R^{*}_{\Pi(t),k}=\sum_{i=1}^{\infty} {\bf I}(X_i(\Pi(t))\geq k), \ \  R_{\Pi(t),k}
=\sum_{i=1}^{\infty} {\bf I}(X_i(\Pi(t))= k).
\]

Let designate $\alpha(x)=\max\{j|\ p_j\geq 1/x\}$.

We assume, following Karlin \cite{21}, that 
\[
\alpha(x)=x^{\theta} L(x), \ \ 0 <\theta< 1.
\]

Here
$L(x)$ is a 
slowly varying function as $ x \to \infty $.
Let for $t\in [0,1], \ k\geq 1$ 
\[
Y_{n,k}^{*}(t) =\frac{R^{*}_{[nt],k}-{\bf E} R^{*}_{[nt],k}}{(\alpha(n))^{1/2}}, 
\ \ \ \ \ \ \ \ \ \  Z_{n,k}^{*}(t) =\frac{R^{*}_{\Pi (nt),k}-{\bf E} R^{*}_{\Pi (nt),k}}{(\alpha(n))^{1/2}},
\]
\[
 Y_{n,k}(t) =\frac{R_{[nt],k}-{\bf E} R_{[nt],k}}{(\alpha(n))^{1/2}},   \ \ \ \ \ \ \ \ \
Z_{n,k}^{**}(t)=\frac{R^{*}_{\Pi ([nt]),k}-{\bf E} R^{*}_{\Pi ([nt]),k}}{(\alpha(n))^{1/2}},
\]
 $K_{k,\theta}=\theta \Gamma(k-\theta)$ for $k>0$ and $ K_{k,\theta}=-\Gamma(1-\theta)$ for $k=0$.

\begin{proposition} \ \ (Theorem 4 in \cite{21}).
Let \ $\theta\in(0,1]$. \ Then $ (R_n-{\bf E} R_n)/B_n^{1/2}$
converges weakly to standard normal distribution, where
\[
B_n=\left\{
\begin{array}{cr}
\Gamma(1-\theta)(2^{\theta}-1)n^{\theta}L(n), & \theta\in (0,1);\\
n \int\limits_0^{\infty}\frac{e^{-1/y}}{y}L(n y)dy \stackrel{def}{=}n L^{*}(n), & \theta=1;
\end{array}
\right.
\]
\end{proposition}

Karlin (\cite{21}, Lemma 4) proved function $L^{*}(x)$ to be slowly varying as $x \to \infty$.

Dutko \cite{14} generalized the theorem by proving asymptotic normality of
$ R_n $ if ${\bf Var} R_n \to \infty$
as $n \to \infty$.  
This condition always holds if $\theta\in(0,1]$ but can hold too for $\theta=0$.
Gnedin, Hansen and Pitman \cite{15} focus on the study of conditions for convergence 
${\bf Var} R_n \to \infty$. They also collect facts on the issue.

\begin{proposition} \  (Theorem 5 in \cite{21}). 
Let $\theta\in(0,1)$, $r_1<\ldots <r_{\nu}$ 
be $\nu$ positive integers. Then random vector
$
\left(Y_{n,r_1}(1), 
\ldots,
Y_{n,r_\nu}(1) \right)
$
converges weakly to a multivariate normal distribution with zero expectation
 and covariances
\[
c_{r_i,r_j}=\left\{
\begin{array}{cr}
-\frac{\theta \Gamma(r_i+r_j-\theta)}{r_i! r_j!} 2^{\theta-r_i-r_j}, & i \neq j;\\
\frac{\theta}{\Gamma(r_i+1)}\left(
\Gamma(r_i-\theta)- 2^{-2r_i+\theta}\frac{\Gamma(2r_i-\theta)}{\Gamma(r_i+1)}\right), & i=j.
\end{array}
\right.
 \]
\end{proposition}

Barbour and Gnedin \cite{5} extended this result to the case $ \theta = 0 $ if
variances go to infinity. They  found conditions for convergence of
covariances to a limit
and identified four types of limiting behavior of variances.
Barbour \cite{4} proved theorems on the approximation of the number of cells with $ k $ balls
by translated Poisson distribution.
Key \cite{22}, \cite{23}  
studied the limit behavior of statistics $R_{n,1}$.
Hwang and Janson \cite{19} 
proved local limit theorems for finite and infinite numbers of cells.
Chebunin \cite{9} constructed $ R_n $-based explicit parameter estimators for a wide range of one-parameter 
families and proved their consistency.

Durieu and Wang \cite{13} established a functional central limit theorem for randomization of a process $R_n$:
indicators are multiplied randomly by $\pm 1$ before summing. The limiting Gaussian process is a sum of independent 
self-similar processes in this case.

Chebunin $\&$ Kovalevskii \cite{10} proved the next FCLT.

\begin{proposition}

{\bf (i)} 
Let $\theta\in(0, 1)$, $\nu\geq 1$ is integer. 
Then process
$
\left(Y^{*}_{n,1}(t),\ldots,Y^{*}_{n,\nu}(t),
\ 0 \leq t \leq 1 \right)
$
converges weakly in the uniform metrics in  $D(0,1)$ 
to $\nu$-dimensional Gaussian process with zero expectation and covariance function 
$(c^{*}_{ij}(\tau,t))_{i,j=1}^{\nu}$: for $\tau\leq t$, $i,j \in \{ 1, \ldots, \nu\}$ (taking $0^0=1$)
\[
c^{*}_{ij}(\tau,t)=\left\{
\begin{array}{cr}
 \sum\limits_{s=0}^{i-1} 
\sum\limits_{m=0}^{j-s-1} \frac{\tau^s (t-\tau)^m K_{m+s,\theta}}{t^{m+s-\theta}s!m!}  
-
 \sum\limits_{s=0}^{i-1} 
\sum\limits_{m=0}^{j-1} \frac{\tau^s t^m K_{m+s,\theta}}{(t+\tau)^{m+s-\theta} s!m!}  , & i < j;\\
t^{\theta} \sum\limits_{m=0}^{j-1} \frac{K_{m,\theta}}{m!}  
-
 \sum\limits_{s=0}^{i-1} 
\sum\limits_{m=0}^{j-1} \frac{\tau^s t^m K_{m+s,\theta}}{(t+\tau)^{m+s-\theta} s!m!}, & i\ge j;\\
\end{array}
\right.
 \]
  $c^*_{ij}(\tau,t)=c^*_{ji}(t,\tau)$.

{\bf (ii)}
Let $\theta=1$. Then process 
$\left( \frac{R_{[nt]}-{\bf E} R_{[nt]}}{(nL^{*}(n))^{1/2}}, \ 0\le t \le 1
\right)$ 
converges weakly in the uniform metrics in $D(0,1)$ 
to a standard Wiener process.
\end{proposition}

\section{Results for Simon model}

Yule showed that, for the Yule---Simon model, 
\begin{equation}\label{Yule_formula}
{\bf E} R_{n,i}/{\bf E} R_n \to f(i), \ \ i\ge 1,
\end{equation}
\[
f(i)=\rho B(i,1+\rho),
\]
$\rho=(1-p)^{-1}$, $B(\cdot,\cdot)$ is Beta function.

We analyze stochastic aspects of this convergence. The limiting distribution is named 
Yule-Simon distribution.

There are many ramifications and applications of the Yule-Simon model. Haight $\&$ Jones \cite{16}
gave special references to word associations tests.
Lansky $\&$ Radill-Weiss \cite{24} proposed a generalization of the model for better
correspondence to applications.

This model can be embedded in a more general context of random cutting 
of recursive trees. In this context, statistics under study are the most frequent words.
See Aldous $\&$ Pitman \cite{2} for its limiting distribution and convergence, 
Baur $\&$ Bertoin \cite{6}, \cite{7} for an overview and new results.
 
Aldous \cite{1} proposed a generalization of the limiting distribution but without 
an underlying process. 

Janson \cite{20} considered generalized Polya urns and proved SLLN, CLT and FCLT for it.
Finite-dimensional vectors 
\[
(R_n, R_{n,1},\ldots, R_{n,m-1}, n - \sum_{i<m} i R_{n,i})
\]
can be studied using these models. So we have componentwise SLLN and finite-dimensional CLT and FCLT 
for these statistics. 

\begin{theorem} For any $p\in (0,1)$ in Simon model and any  $m> 1$ 
\[
\frac{R_n}{n} \to p \ \ {\rm a.s.},
\]
\[
\frac{(R_{n,1}, \ldots, R_{n,m-1})}{n} \to \frac{p}{1-p}\left(B\left(1, \frac{2-p}{1-p}\right),\ldots 
B\left(m-1, \frac{2-p}{1-p}\right)\right) \ \ {\rm a.s.},
\]
and, in $D[0,\infty)$,
\[
\left\{
n^{-1/2}(
R_{[nt]}-tnp), \right.
\]
\[
\left.  
n^{-1/2}\left(R_{[nt],j} - t n \frac{p}{1-p} B\left(j, \frac{2-p}{1-p}\right)\right), \ 1 \le j \le m-1, \ t \ge 0 
\right\} \to_{d} V(t),
\]
where $V(t)$ is a continuous centered $m$-dimensional Gaussian process,
its covariance matrix-function ${\bf E} V(x)V^T(y)$ depends on $p,x,y$ only.
\end{theorem}

{\it Proof}

Simon's model can be studied as a very partial case of Janson's \cite{20} urn scheme. 
In this model, Simon urns with $1\le i\le m-1$ balls are balls with numbers from 2 to $m$ with weights $a_i=i-1$.
The $(m+1)$-th urn contains all other balls with weights $a_{m+1}=1$. 
Balls with number 1 correspond to all different words, $a_1=0$.
So the random uniform choice of balls in the Simon model corresponds to the random choice with weights $a_i$, 
$1\le i \le m+1$, in Janson model.

Let $({\bf e}_1,\ldots,{\bf e}_{m+1})$ be the standard basis in ${\bf R}^m$, then increment vectors in Janson model are
\[
\xi_i={\bf e}_1+{\bf e}_2 \ {\rm with \ probability} \ p, \ \   \xi_i={\bf e}_{i+1}-{\bf e}_{i} \ {\rm with \ probability} \ q=1-p,
\ \ 2\le i<m,
\]
\[
\xi_{m}={\bf e}_1+{\bf e}_2  \ {\rm with \ probability} \ p, \ \   \xi_{m}=m{\bf e}_{m+1}-{\bf e}_{m} \ {\rm with \ probability} \ q=1-p,
\]
\[
\xi_{m+1}={\bf e}_1+{\bf e}_2  \ {\rm with \ probability} \ p, \ \   \xi_{m+1}={\bf e}_{m+1} \ {\rm with \ probability} \ q=1-p.
\]

So matrix $A=(a_j {\bf E} \xi_{ji})$ has the first (maximal) eigenvalue $\lambda_1=1$, the next (second)
eigenvalue $\lambda_2=0$. Thus the assumption $\lambda_2<\lambda_1/2$ of Theorem 2.22 and Theorem 2.31(i) 
in Janson \cite{20} is hold.
So we use Janson's theorems 2.21 (SLLN), 2.22, 2.31(i).

The proof is complete.

\section{Offline Simon modification of the Karlin model}

Let the balls take urns as in the Karlin model, but the appearance of a new non-empty urn corresponds only to
the embryo of a new word.

After the finish of the Karlin model, Simon's model starts to work in moments of the appearance of the embryos of new words:
the $ k $-th embryo coincides with each of the previous embryos with probability $ q / k $, and is a new
word with probability $ p $, $ p + q = 1 $. 
We use Janson's theorems.

Balls of the first type are embryos of words encountered exactly once, their number at the $ k $ -th step is $ X_{k, 1} $, balls
of this type have weights (responsible for the probabilities of choosing balls) $ a_1 = 1 $.

Balls of the second type are embryos of words encountered more than once, their number at the $ k $ -th step is $ X_{k, 2} $, balls
of this type  have weights $ a_2 = 1 $ also.

The system starts from the state $ X_{0,1} = 1, X_{0,2} = 0 $. At each step $ k $, one ball  is chosen
at random from $ X_ {k, 1} + X_ {k, 2} $ balls of the first and second types.

If the selected ball is of the first type, then a ball of the first type is added with probability $ p $,
and it is replaced by a ball of the second type with probability~$ q $:
\[
X_0=(1,0)',
\]
\[
\xi_1=(1,0)' \ \text{with probability} \ p ,
\]
\[
\xi_1=(-1,1)' \ \text{with probability} \ q .
\]

If the chosen ball is of the second type, then a ball of the first type is added with probability $ p $,
and  nothing happens with probability $ q $:
\[
\xi_2=(1,0)' \ \text{with probability} \ p ,
\]
\[
\xi_2=(0,0)' \ \text{with probability} \ q .
\]

So
\[
A=\left(a_j {\bf E} \xi_{ji}\right)_{i,j=1}^2=
\left(
\begin{array}{cc}
p-q & p \\
q & 0 
\end{array}
\right).
\]

Eigenvalues of $A$ are $\lambda_1=p$, $\lambda_2=-q<\lambda_1/2$.

Eigenvector  $v_1$ of $A$ correspond to eigenvalue $\lambda_1=p$ 
and condition (2.2) from Janson \cite{20} $a\cdot v_1=1$.
As $a=(1 \, 1 )'$, we have
\[
v_1=\left(
\begin{array}{c}
p\\
q
\end{array}
\right).
\]

We calculate $e^{tA}$ using Sylvester's formula:
\[
e^{tA}=\frac{e^{pt}}{p+q}(A+qI)- \frac{e^{-qt}}{p+q}(A-pI)
=e^{pt}(A+qI)- e^{-qt}(A-pI).
\]

From (3.17) in Janson \cite{20},
\[
\phi(s,A)=\int_0^s e^{tA} dt 
=\frac{e^{ps}-1}{p}(A+qI)+ \frac{e^{-qs}-1}{q}(A-pI)
\]

From (3.18) in \cite{20},
\[
\psi(s,A)=e^{sA}-\lambda_1 v_1 a' \phi(s,A)
\]
\[
=e^{ps}(A+qI)- e^{-qs}(A-pI)
- p v_1 a' \left(\frac{e^{ps}-1}{p}(A+qI)+ \frac{e^{-qs}-1}{q}(A-pI)\right).
\]

Note that
\[
v_1 a'=A+qI=(A+qI)^2.
\]

So
\[
\psi(s,A)
=- e^{-qs}(A-pI)
+A+qI- \frac{p(e^{-qs}-1)}{q}(A+qI)(A-pI).
\]
\[
=A+qI- e^{-qs}(A-pI).
\]

Let
\[
A_0:=A+qI,
\]
then
\[
\psi(s,A)
=A_0(1- e^{-qs})+ Ie^{-qs}.
\]

We calculate matrices $B_1$, $B_2$, $B$. From (2.13) in \cite{20}, 
\[
B_1={\bf E} \xi_1 \xi_1' = 
p 
\left(
\begin{array}{cc}
1 & 0  \\
0 & 0 
\end{array}
\right)+
q
\left(
\begin{array}{cc}
1 & -1  \\
-1 & 1 
\end{array}
\right)
=
\left(
\begin{array}{cc}
1 & -q \\
-q & q
\end{array}
\right),
\]
\[
B_2={\bf E} \xi_2 \xi_2' = 
\left(
\begin{array}{cc}
p & 0  \\
0 & 0 
\end{array}
\right).
\]

From (2.14) in \cite{20},
\[
B=a_1 v_{11} B_1+a_2 v_{12} B_2=p
\left(
\begin{array}{cc}
1+q & -q \\
-q & q
\end{array}
\right).
\]

From Theorem 3.21 in \cite{20},
\[
\frac{X_n}{n} \to_{a.s.} \lambda_1 v_1=p v_1.
\]

From Theorem 3.22 in \cite{20},
\[
\frac{X_n-n \lambda_1 v_1}{\sqrt{n}} \to_{d} N(0,\Sigma),
\]
\[
\Sigma=\int_0^{\infty} \psi(s,A)B\psi(s,A)'e^{-\lambda_1 s}\lambda_1 ds - \lambda_1^2 v_1 v_1'
\]
\[
=\int_0^{\infty} \left(A_0(1- e^{-qs})+ Ie^{-qs}\right)B\left(A_0'(1- e^{-qs})+ Ie^{-qs}\right) e^{-ps}p ds
\]
\[
-p^2
\left(
\begin{array}{cc}
p^2 & pq \\
pq & q^2
\end{array}
\right)
\]
\[
=\int_0^{\infty} \left(A_0 B A_0' + e^{-qs}(A_0B+BA_0'-2A_0 B A_0')+e^{-2qs}(B+A_0 B A_0')
\right) e^{-ps}p ds
\]
\[
-p^2
\left(
\begin{array}{cc}
p^2 & pq \\
pq & q^2
\end{array}
\right)
\]
\[
=A_0 B A_0' +p(A_0B+BA_0'-2A_0 B A_0')+\frac{p}{1+q}(B+A_0 B A_0')
-p^2
\left(
\begin{array}{cc}
p^2 & pq \\
pq & q^2
\end{array}
\right)
\]
\[
=pq
\left(
\begin{array}{cc}
3p^2+p & p-3p^2 \\
p-3p^2 & 1-3p+3p^2
\end{array}
\right).
\]

From Theorem 3.31 in \cite{20}, there is FCLT for $\{X_{[nt]}, \ 0 \le t \le 1\}$:
\[
\left\{ \frac{X_{[nt]}-nt \lambda_1 v_1}{\sqrt{n}}, \ 0 \le t \le 1 \right\} \to_{d} U=\{U(t), \ 0 \le t \le 1\}
\]
 in the uniform metrics, $U$ is the centered Gaussian process, and for $0<x\le y$
 \[
 {\bf E} U(x)U'(y)=\int\limits_{-\lambda_1^{-1}\log x}^{\infty} \psi(s+\lambda_1^{-1}\log x, A)B\psi(s+\lambda_1^{-1}\log y, A)' e^{-\lambda_1 x} \lambda_1 dx - x \lambda_1^2 v_1 v_1'
 \]
 \[
=\int\limits_{-\lambda_1^{-1}\log x}^{\infty} 
\left(A_0(1- e^{-q(s+\log x/p)})+ Ie^{-q(s+\log x/p)}\right) 
B
\]
\[
\times
\left(A_0'(1- e^{-q(s+\log y/p)})+ Ie^{-q(s+\log y/p)}\right) 
e^{-\lambda_1 x} \lambda_1 dx - x \lambda_1^2 v_1 v_1'
 \]
\[
=xpq
\left(
\begin{array}{cc}
2p^2 & 2pq \\
pq-p^2 & q^2-pq
\end{array}
\right)
+
xp^2q \left(\frac{x}{y}\right)^{q/p}
\left(
\begin{array}{cc}
1+p & -1-p \\
-p & p
\end{array}
\right).
\]

We have a good match for {\it Childe Harold's Pilgrimage} with $\theta=0.58$, 
$p=0.56$,
for {\it Evgene Onegin} with $\theta=0.74$, 
$p=0.71$.

\section{Compound Poisson modification of the Karlin model}

We want to invent a text model such that $\theta^*$ converges to a number less than $\widehat{\theta}$.
We propose the following simple model: words appear in the same way as in Karlin's model, but
 each word with probability $q_i$ (independently of others and of the process of the appearance of words) is repeated $i\ge 1$ times, $q_1+q_2+\ldots=1$. New statistics of the number of different words and the number of words that occur once are denoted by
$\mathbb{R}_n$ and $\mathbb{R}_{n,1}$. 

The Poisson version of these processes uses the compound Poisson process.

We want to invent a text model such that $\theta^*$ converges to a number less than $\widehat{\theta}$.
We propose the following simple model: words appear in the same way as in Karlin's model, but
 each word with probability $q_i$ (independently of others and of the process of the appearance of words) is repeated $i\ge 1$ times, $q_1+q_2+\ldots=1$. New statistics of the number of different words and the number of words that occur once are denoted by ${\mathbb R}_n$ and ${\mathbb R}_{n,1}$. We denote by $ {\mathbb X}_i (n) $ a number of balls in urn $ i $, and
$$
\mathbb{R}_{n, k}^{*}=\sum_{i=1}^{\infty} \mathbf{I}\left(\mathbb{X}_i(n)\ge k\right), \ \ \mathbb{R}_{n}=\mathbb{R}_{n, 1}^{*}, \ \ 
\mathbb{R}_{n, k}=\mathbb{R}_{n, k}^{*}-\mathbb{R}_{n, k+1}^{*}.
$$

The Poisson version of these processes uses the compound Poisson process. Let $ \Pi=\{\Pi(t), t \geq 0\}$ be a Poisson process with parameter $1$. The Poissonized version of Karlin model assumes the total number of $\Pi(n)$ balls. According to well-known thinning property of Poisson flows, stochastic processes $\left\{{\mathbb X}_{i}(\Pi(t)) \stackrel{\text { def }}{=} \mathbf{\Pi}_{i}(t), t \geq 0\right\}$ are compound Poisson with intensities $p_{i}$ and are mutually independent for different $i$ 's. The definition implies that
$$
\mathbb{R}_{\Pi(n), k}^{*}=\sum_{i=1}^{\infty} \mathbf{I}\left(\mathbf{\Pi}_{i}(n) \geq k\right), \quad \mathbb{R}_{\Pi(n), k}=\sum_{i=1}^{\infty} \mathbf{I}\left(\mathbf{\Pi}_{i}(n)=k\right) .
$$

We proved FCLT for the vector process $({\mathbb R}_n, \ {\mathbb R}_{n,1})$. Let $\alpha(x)=\max \left\{j \mid p_{j} \geq 1 / x\right\} .$ Following Karlin \cite{21}, we assume that $\alpha(x)=x^{\theta} L(x), 0 \leq \theta \leq 1$. Here $L(x)$ is a slowly varying function as $x \rightarrow \infty$. Let for $t \in[0,1], k \geq 1$
\[
\mathbb{Y}_{n, k}^{*}(t)=\frac{\mathbb{R}_{[n t], k}^{*}-\mathbf{E} \mathbb{R}_{[n t], k}^{*}}{(\alpha(n))^{1 / 2}}, \ \  \mathbb{Z}_{n, k}^{*}(t)=\frac{\mathbb{R}_{\Pi(n t), k}^{*}-\mathbf{E} \mathbb{R}_{\Pi(n t), k}^{*}}{(\alpha(n))^{1 / 2}},
\]
\[
\ \  \mathbb{Y}_{n, k}(t)=\frac{\mathbb{R}_{[n t], k}-\mathbf{E} \mathbb{R}_{[n t], k}}{(\alpha(n))^{1 / 2}}.
\]

\begin{theorem}
  Let $\theta \in(0,1)$ and $v \geq 1$ be an integer. Then process $\left(\mathbb{Y}_{n, 1}^{*}(t), \ldots, \mathbb{Y}_{n, v}^{*}(t), 0 \leq t \leq 1\right)$ converges weakly in the uniform metric in $D\left([0,1]^{v}\right)$ to $v$ -dimensional Gaussian process with zero expectation and covariance function $\left(\mathbb{C}_{i j}^{*}(\tau, t)\right)_{i, j=1}^{v}$. 
\end{theorem}

\begin{lemma}
{\bf (i)} There exist $n_0\ge 1$, $C(\theta)<\infty$ such that 
$
{\bf E} R_{\Pi(n\delta)}/\alpha(n) 
\le
C(\theta) \delta^{\theta/2}
$
for any $\delta\in[0,1], \ n\ge n_0$.

{\bf (ii)}
Let $\tau\le t$, then ${\bf E} (\mathbb{R}^*_{\Pi(t),k}-\mathbb{R}^*_{\Pi(\tau),k})\le {\bf E} R_{\Pi(t-\tau)}$, $k\ge1$.

{\bf (iii)}
For any pair $\varepsilon$,
  $\delta\in(0,1)$ there exists integer $n_0$ such that 
$
{\bf P}(\forall t\in[0,1] \  \ \exists \tau: |\tau-t|\le \delta,  
\ \Pi(n\tau)= [nt]) \stackrel{def}{=} {\bf P} (A(n))\ge1-\varepsilon/2
$
for any $n\ge n_0$. \rm
\end{lemma}

{\it Proof.} {\bf (i)} and  {\bf (iii)} were proved in Chebunin and Kovalevskii \cite{12},
{\bf (ii)} follows from
\[
\mathbf{E}\left(\mathbb{R}_{\Pi(t), k}^{*}-\mathbb{R}_{\Pi(\tau), k}^{*}\right)=\sum_{i=1}^{\infty} \sum_{j=0}^{k-1} \mathbf{P}\left(\mathbf{\Pi}_{i}(\tau)=j\right) \mathbf{P}\left(\mathbf{\Pi}_{i}(t)-\mathbf{\Pi}_{i}(\tau) \geq k-j\right) 
\]
\[
\leq \sum_{i=1}^{\infty} \mathbf{P}\left(\mathbf{\Pi}_{i}(t-\tau) \geq 1\right)
=\mathbf{E} \mathbb{R}_{\Pi(t-\tau)}= \mathbf{E} R_{\Pi(t-\tau)}.
\]

The proof is complete.

{\it Proof of Theorem~2}

{\bf Step 1 (covariances)} Let $\tau \leq t$, $i\le j$ and $\pi_{k,i}(t)=\mathbf{P}(\mathbf{\Pi}_k(t)=i)$
\[
{\bf cov}\left(\mathbb{R}_{\Pi(\tau), i}^{*}, \mathbb{R}_{\Pi(t), j}^{*}\right) 
=\sum_{k=1}^{\infty}\left(\mathbf{P}\left(\mathbf{\Pi}_{k}(\tau)<i, \mathbf{\Pi}_{k}(t)<j\right)-\mathbf{P}\left(\mathbf{\Pi}_{k}(\tau)<i\right) \mathbf{P}\left(\mathbf{\Pi}_{k}(t)<j\right)\right)
\]
\[
=\sum_{k=1}^{\infty}\sum_{s=0}^{i-1}\pi_{k,s}(\tau)\left(\sum_{m=0}^{j-s-1}\pi_{k,m}(t-\tau)-\sum_{m=0}^{j-1}\pi_{k,m}(t)\right),
\]
and $i>j$
\[
{\bf cov}\left(\mathbb{R}_{\Pi(\tau), i}^{*}, \mathbb{R}_{\Pi(t), j}^{*}\right) 
=\sum_{k=1}^{\infty}\sum_{s=0}^{j-1}\pi_{k,s}(t)\left(1-\sum_{m=0}^{i-1}\pi_{k,m}(\tau)\right).
\]
Let $t'= t(1-q_0)$, $c_1=q_1/(1-q_0)$
\[
\mathbf{E} \mathbb{R}_{\Pi(t)}=\sum_{k=1}^{\infty}1-\pi_{k,0}(t)=
\sum_{k=1}^{\infty}1-e^{-p_k t}=\mathbf{E}R_{\Pi(t)}, \ \ 
\mathbf{D} \mathbb{R}_{\Pi(t)}=\mathbf{D}R_{\Pi(t)}.
\]
\[
\mathbf{E} \mathbb{R}_{\Pi(t),1}=\sum_{k=1}^{\infty}\pi_{k,1}(t)=
\sum_{k=1}^{\infty} q_1\mathbf{P}(\Pi_k(t)=1) 
=q_1\mathbf{E}R_{\Pi(t) ,1}.
\]

\[
{\bf cov}\left(\mathbb{R}_{\Pi(\tau)}^{*}, \mathbb{R}_{\Pi(t)}^{*}\right) 
=\sum_{k=1}^{\infty}\pi_{i,0}(t)(1-\pi_{i,0}(\tau))=\mathbf{E} R_{\Pi(t+\tau)}-\mathbf{E} R_{\Pi(t)}.
\]
\[
{\bf cov}\left(\mathbb{R}_{\Pi(\tau), 2}^{*}, \mathbb{R}_{\Pi(t), 2}^{*}\right) 
=\sum_{k=1}^{\infty}\sum_{s=0}^{1}\pi_{k,s}(\tau)\left(\sum_{m=0}^{1-s}\pi_{k,m}(t-\tau)-\sum_{m=0}^{1}\pi_{k,m}(t)\right)
\]
\[
=\sum_{k=1}^{\infty}\pi_{i,0}(t)(1-\pi_{i,0}(\tau))+e^{-p_k \tau}(q_1p_k(t-\tau)e^{-p_k(t-\tau)}-q_1p_kte^{-p_kt})
\]
\[
+\sum_{k=1}^{\infty}q_1p_k\tau e^{-p_k \tau}(e^{-p_k(t-\tau)}-e^{-p_kt}-q_1p_kt e^{-p_kt})
\]
\[
=\mathbf{E} (R_{\Pi(t+\tau)}- R_{\Pi(t)})+q_1\mathbf{E}R_{\Pi(t) ,1}
-q_1 \mathbf{E}R_{\Pi(t+\tau) ,1}-\frac{2q_1^2 t\tau}{(t+\tau)^2}\mathbf{E}R_{\Pi(t+\tau) ,2}.
\]
\[
{\bf cov}\left(\mathbb{R}_{\Pi(\tau)}^{*}, \mathbb{R}_{\Pi(t),2}^{*}\right) 
=\sum_{k=1}^{\infty}\pi_{k,0}(\tau)\left(\sum_{m=0}^{1}\pi_{k,m}(t-\tau)-\sum_{m=0}^{1}\pi_{k,m}(t)\right)
\]
\[
=\mathbf{E} (R_{\Pi(t+\tau)}- R_{\Pi(t)})+\frac{q_1(t-\tau)}{t}\mathbf{E}R_{\Pi(t) ,1}
-\frac{q_1 t}{t+\tau} \mathbf{E}R_{\Pi(t+\tau) ,1}.
\]
\[
{\bf cov}\left(\mathbb{R}_{\Pi(t)}^{*}, \mathbb{R}_{\Pi(\tau),2}^{*}\right) 
=\sum_{k=1}^{\infty}\pi_{k,0}(t)\left(1-\pi_{k,0}(\tau)-\pi_{k,1}(\tau)\right)
\]
\[
=\mathbf{E} (R_{\Pi(t+\tau)}- R_{\Pi(t)})
-\frac{q_1 \tau}{t+\tau} \mathbf{E}R_{\Pi(t+\tau) ,1}.
\]
Since 
 \[
{\bf E} R_{\Pi (t)}\sim \Gamma(1-\theta) \alpha(t), \ {\bf E} R_{\Pi (t),k}\sim \theta\frac{\Gamma(k-\theta)}{k!} \alpha(t)\ \ \textrm{if} \ \ \theta\in(0,1), \ k\ge1,
\]
then
\[
\mathbb{C}_{1 1}^{*}(\tau, t)=\Gamma(1-\theta)((t+\tau)^\theta-t^\theta), 
\]
\[
\ \mathbb{C}_{2 2}^{*}(\tau, t)=\Gamma(1-\theta)\left(((t+\tau)^\theta-t^\theta)(1-q_1\theta) -\frac{q_1^2\theta(1-\theta) t\tau}{(t+\tau)^{2-\theta}}\right), 
\]
\[
\mathbb{C}_{1 2}^{*}(\tau, t)=\Gamma(1-\theta)\left((t+\tau)^\theta-t^\theta + \frac{q_1\theta (t-\tau)}{t^{1-\theta}}-\frac{q_1 \theta t}{(t+\tau)^{1-\theta}}\right),
\]
\[
\mathbb{C}_{1 2}^{*}(t, \tau)=\Gamma(1-\theta)\left((t+\tau)^\theta-t^\theta -\frac{q_1 \theta \tau}{(t+\tau)^{1-\theta}}\right).
\]

{\bf  Step~2 (convergence of finite-dimensional distributions)}
Analogously to proof of Theorem 1 in \cite{14}, we have that, for any fixed $m \geq 1,0<t_{1}<t_{2}<\cdots<t_{m} \leq 1$ triangle array of $m v$ -dimensional random vectors 
$\left\{
( ( \mathbf{I}(\mathbf{\Pi}_k(n t_{j}) \geq i)-\mathbf{P}(\mathbf{\Pi}_{k}(n t_{j}) \geq i)) \alpha^{-1 / 2}(n), i \leq v, j \leq m), k \leq n
\right\}_{n \geq 1}$ 
satisfies Lindeberg condition.

{\bf Step~3 (relative compactness)}
 Let for any $\tau_{1} \leq \tau_{2}$,
\[
\mathbb{R}_{\Pi\left(\tau_{2}\right), k}^{*}-\mathbb{R}_{\Pi\left(\tau_{1}\right), k}^{*}=\sum_{i=1}^{\infty} \mathbf{I}\left(\mathbf{\Pi}_{i}\left(\tau_{2}\right) \geq k, \mathbf{\Pi}_{i}\left(\tau_{1}\right)<k\right) \stackrel{\text { def }}{=} \sum_{i=1}^{\infty} \mathbb{I}_{i}\left(\tau_{1}, \tau_{2}\right)=\sum_{i=1}^{\infty} \mathbb{I}_{i},
\]
$\mathbb{P}_{i}=\mathbb{P}_{i}\left(\tau_{1}, \tau_{2}\right)=\mathbf{P}\left(\mathbb{I}_{i}\right) .$ We will use designations $\mathbb{I}_{i}$ and corresponding $\mathbb{P}_{i}$ for different values of $\tau_{1}<\tau_{2}$.

We need a new process $\mathbb{Z}_{n, k}^{* *}(t)=\frac{\mathbb{R}_{\Pi([n t]), k}^{*}-\mathbf{E} \mathbb{R}_{\Pi([n t]), k}^{*}}{(\alpha(n))^{1 / 2}}$.

We (a) prove continuity of the limiting process; (b) prove that $\mathbb{Z}_{n, k}^{*}$ and $\mathbb{Z}_{n, k}^{* *}$ are 'close'; (c) prove relative compactness of $\mathbb{Z}_{n, k}^{* *}$.

(a) Let $\tau_{1}=n t_{1}, \tau_{2}=n t_{2}$ for $t_{1}<t_{2},$ then
\[
\mathbf{E}\left(\mathbb{Z}_{n, k}^{*}\left(t_{2}\right)-\mathbb{Z}_{n, k}^{*}\left(t_{1}\right)\right)^{2}
=\sum_{i=1}^{\infty} \mathbf{E}\left(\mathbb{I}_{i}-\mathbb{P}_{i}\right)^{2} / \alpha(n) 
\]
\[
\leq \sum_{i=1}^{\infty} \mathbb{P}_{i} / \alpha(n) \leq C(\theta)\left(t_{2}-t_{1}\right)^{\theta / 2}.
\]
Here we used the fact that variance of an indicator is lesser than its expectation and Lemma 1 (i, ii). Using Step 1 and Theorem 1.4 in Adler \cite{0}, we prove that the $k$-th component of the limiting Gaussian process is in $C(0,1)$ a.s. As the limiting Gaussian process is in $C\left([0,1]^{v}\right)$ a.s., weak convergence in Skorokhod topology implies the same in the uniform topology.

(b) As $\mathbb{R}_{\Pi(n t), k}^{*}-\mathbb{R}_{\Pi([n t]), k}^{*} \leq \Pi([n t]+1)-\Pi([n t])$ a.s., and $\mathbf{E}(\Pi([n t]+1)-\Pi([n t]))=1$ we have for any $\eta>0$
\[
\mathbf{P}\left(\sup _{0 \leq t \leq 1}\left|\mathbb{Z}_{n, k}^{*}(t)-\mathbb{Z}_{n, k}^{* *}(t)\right|>\eta\right) 
\]
\[ \leq \mathbf{P}\left(\sup _{0 \leq m \leq n}(\Pi(m+1)-\Pi(m)+1)>\eta \sqrt{\alpha(n)}\right) 
\]
\[
\leq \sum_{m=0}^{n} \mathbf{E} e^{\Pi(m+1)-\Pi(m)+1} / e^{\eta \sqrt{\alpha(n)}}=(n+1) e^{e-\eta \sqrt{\alpha(n)}} \rightarrow 0
\]
as $n \rightarrow \infty$. So it is enough to prove relative compactness of $\left\{\mathbb{Z}_{n . k}^{* *}\right\}_{n \geq n_{0}}$.

(c) Let $t_{2}-t_{1} \geq \frac{1}{2 n},$ then $\left[n t_{2}\right]-\left[n t_{1}\right] \leq n\left(t_{2}-t_{1}\right)+1 \leq 3 n\left(t_{2}-t_{1}\right) .$ Let $\gamma=[16 / \theta]+1,$ and $\tau_{1}=\left[n t_{1}\right], \tau_{2}=\left[n t_{2}\right]$.
Using independence of terms and Rosenthal inequality, we have for all $n \geq n_{0}$ (where $n_{0}$ is from Lemma 1 (i))
\[
\mathbf{E}\left|\mathbb{Z}_{n, k}^{* *}\left(t_{2}\right)-\mathbb{Z}_{n, k}^{* *}\left(t_{1}\right)\right|^{\gamma}  \leq \frac{c(\gamma)}{(\alpha(n))^{\gamma / 2}}\left(\sum_{i=1}^{\infty} \mathbf{E}\left|\mathbb{I}_{i}-\mathbb{P}_{i}\right|^{\gamma}+\left(\sum_{i=1}^{\infty} \mathbf{E}\left(\mathbb{I}_{i}-\mathbb{P}_{i}\right)^{2}\right)^{\gamma / 2}\right) 
\]
\[
 \leq \frac{c(\gamma)}{(\alpha(n))^{\gamma / 2}}\left(\sum_{i=1}^{\infty} \mathbb{P}_{i}+\left(\sum_{i=1}^{\infty} \mathbb{P}_{i}\right)^{\gamma / 2}\right) 
 \]
 \[
\leq \frac{c(\gamma)}{(\alpha(n))^{\gamma / 2}}\left(24 n^{4}\left(t_{2}-t_{1}\right)^{4}+\left(\mathbf{E} R_{\Pi\left(3 n\left(t_{2}-t_{1}\right)\right)}\right)^{\gamma / 2}\right) \leq \widetilde{C}(\theta)\left(t_{2}-t_{1}\right)^{4}.
\]

Here $c(\gamma)$ and $\widetilde{C}(\theta)$ depend on its argument only. Above we used the fact that variance of an indicator is lesser than its expectation, inequality $\sum_{i} \mathbb{P}_{i} \leq \mathbf{E}\left(\Pi\left(\left[n t_{2}\right]\right)-\Pi\left(\left[n t_{1}\right]\right)\right) \leq 3 n\left(t_{2}-t_{1}\right) \leq 24 n^{4}\left(t_{2}-t_{1}\right)^{4},$ and Lemma 1(i, ii).

Let $0 \leq t_{2}-t_{1}<1 / n,$ then $\left[n t_{1}\right]=[n t]$ or $\left[n t_{2}\right]=[n t]$ for any $t \in\left[t_{1}, t_{2}\right] .$ So $\mathbb{Q} \stackrel{\text { def }}{=} \mathbf{E}\left(\left|\mathbb{Z}_{n, k}^{* *}(t)-\mathbb{Z}_{n, k}^{* *}\left(t_{1}\right)\right|^{\gamma / 2} \mid \mathbb{Z}_{n, k}^{* *}\left(t_{2}\right)-\right.$
$\left.\left.\mathbb{Z}_{n, k}^{* *}(t)\right|^{\gamma / 2}\right)=0 \leq\left(t_{2}-t_{1}\right)^{2}$.

Let $t_{2}-t_{1} \geq 1 / n,$ then there are 3 possible cases:

(1) $t_{2}-t \geq \frac{1}{2 n}, t-t_{1} \geq \frac{1}{2 n},$ then from Cauchy-Bunyakovsky Inequality, $\mathbb{Q} \leq \widetilde{C}(\theta)\left(t_{2}-t\right)^{2}\left(t-t_{1}\right)^{2} \leq \widetilde{C}(\theta)\left(t_{2}-t_{1}\right)^{2}$.

(2) $t_{2}-t \geq \frac{1}{2 n}, t-t_{1}<\frac{1}{2 n},$ then from Cauchy-Bunyakovsky Inequality,
$$
\mathbb{Q} \leq \sqrt{\widetilde{C}(\theta)\left(t_{2}-t\right)^{4} \mathbf{E}\left(\frac{\Pi(1)+1}{\sqrt{\alpha(n)}}\right)^{\gamma}} \leq \widehat{C}(\theta)\left(t_{2}-t_{1}\right)^{2}.
$$

(3) $t_{2}-t<\frac{1}{2 n}, t-t_{1} \geq \frac{1}{2 n}$, symmetric to case 2. 

So we have (see Billingsley \cite{8+}, Theorem 13.5) tihgtness of $k$-th component and therefore tihgtness of all the vector.

{\bf Step 4 (approximation of the original process) }
From the relative compactness of distributions of processes $\left\{\mathbb{Z}_{n, k}^{*}\right\}_{n \geq n_{0}, k \geq 1}$ we get that for every pair $\varepsilon>0, \eta>0$ there exist $\delta \in(0,1)$ and $N_{1}=N_{1}(\varepsilon, \eta)$ such that for all $n \geq N_{1}$
$$
\mathbf{P}\left(\sup _{|t-\tau| \leq \delta}\left|\mathbb{Z}_{n, k}^{*}(\tau)-\mathbb{Z}_{n, k}^{*}(t)\right| \geq \eta\right) \leq \varepsilon
$$
Then $\left(\right.$ as $\left.\mathbf{P}\left(\mathbb{Y}_{n, k}^{*}(t)=\mathbb{Z}_{n, k}^{*}(\tau) \mid \Pi(n \tau)=[n t]\right)=1\right)$ we have for all $n \geq \max \left(N, N_{1}\right),$ where $N$ is from Lemma 1 (iii),
$$
\begin{aligned}
\mathbf{P}\left(\sup _{0 \leq t \leq 1}\left|\mathbb{Y}_{n, k}^{*}(t)-\mathbb{Z}_{n, k}^{*}(t)\right| \geq \eta\right) & \leq \mathbf{P}\left(\sup _{0 \leq t \leq 1}\left|\mathbb{Y}_{n, k}^{*}(t)-\mathbb{Z}_{n, k}^{*}(t)\right| \geq \eta, A(n)\right)+\varepsilon \\
& \leq \mathbf{P}\left(\sup _{|t-\tau| \leq \delta}\left|\mathbb{Z}_{n, k}^{*}(\tau)-\mathbb{Z}_{n, k}^{*}(t)\right| \geq \eta\right)+\varepsilon \leq 2 \varepsilon.
\end{aligned}
$$

The proof is complete.

For {\it Childe Harold's Pilgrimage}  we have a good match with $\theta=0.54$, 
$q_1=1$.
For {\it Evgene Onegin} we need some nonzero $q_2$, we have $\theta=0.72$, 
$q_1=0.96$, $q_2=0.04$.

\section{Online Simon modification of the Karlin model}

The disadvantage of the Simon model is the linear growth of $R_n$. We need a power growth with an exponent lesser than $1$.
Our idea is to use the classical infinite urn model with re-distribution: any ball takes an urn independently 
with some discrete power law; if the ball is falling in an empty urn then 
 it is re-distributed uniformly on all previous urns with one ball with probability $1-p$, and it stays in the new urn with probability $p$. 

The model starts from the infinite sequence of empty urns. 
The first ball takes one of the urns with the integer-valued power law with exponent $1/\beta$,
$0<\beta<1$. Each next ball takes one of the urns with the same law, independently of previous balls. 
After this, any next ball, if it is in a new urn, is re-tossed independently, that is, with  probability  $1-p $ 
selects one of the previous urns with one ball at random and joins
it, like in the Simon model. 
All other balls stay in selected urns. This is the online variant of the model from Section 5. 

We have a good match for {\it Childe Harold's Pilgrimage} with $\beta=0.54$, $p=0$.
We need some nonzero $p$ to correspond the model to {\it Evgene Onegin}, $\beta=0.73$, $p=0.08$.

The simulation shows that statistics $\widehat{\theta}$ and $\theta^{*}$  converge to different limits under these models, 
in contrast to the elementary urn model.
However,  the analytical dependencies of these limits on the parameters of the models remain unclear.

\section{Discussion}

So, the infinite urn model (Karlin model) 
and the Simon model are the most known models of texts, but they have disadvantages that do not give the ability to 
simulate the number of unique words correctly. The infinite urn model states too strict conditions for the relative number of unique and different words.
The Simon model states a linear growth of the numbers of different and unique words. Its modification by Baur and Bertoin  \cite{8} preserves the 
linear dependence. 

We propose three modifications of the Karlin and Simon models. The first one is the offline variant, the Simon model starts after the completion of the infinite 
urn scheme. We have analytical theorems (SLLN and FCLT) for the process in embedded times of increation of the initial urn process.
This model holds any relative number of unique and different words, but we can test the model under the condition of the given initial urn process only. 

The second modification involves the compound Poisson process in the infinite urn model. As a variant, any word can be doubled with some 
positive probability. We prove SLLN and FCLT for the modification. This model can give the decrease (not increase) of the number of different words 
only. But this decrease is helpful in applications. So we have a simple model that covers the range of parameters that is actual for 
applications. We have analytical results that can be used for testing the model. On the other side, the model with repeating words is strange.

The third modification is the online variant, the Simon redistribution works at any toss of the Karlin model.  Similar to the compound Poisson model, 
it can give the decrease of the number of unique words only.  But, in contrast to the compound Poisson model, we have no analytics
(limit theorems) for this modification.
The online modification seems to be the most logically supported, but we can study it by simulation only due to the absence of theoretically supported tests.

\bigskip

{\bf Acknowledgements}
The reported study was funded by RFBR and CNRS according to the research project No. 19-51-15001.
The authors like to thank 
Bastien Mallein for 
fruitful and constructive comments and suggestions, especially for the interpretation of Simon's model 
in terms of Janson's generalized Polya urn scheme.
The authors thank 
Sergey Foss and Sanjay Ramassamy for helpful discussions.

\end{document}